\documentclass[draft]{amsart}

\usepackage{amssymb}

\usepackage[lowtilde]{url}

\theoremstyle{plain}

\numberwithin{theorem}{section}
\numberwithin{equation}{section}

\numberwithin{enumi}{equation}

\newcommand\thmcall[1]{
  \setcounter{theorem}{\value{equation}}
  \numberwithin{equation}{theorem}    
  \numberwithin{enumi}{theorem}       
  \begin{#1}
  }

\newcommand\exitthmcall[1]{
    \end{#1}
    \numberwithin{equation}{section}  
    \numberwithin{enumi}{equation}    
    \setcounter{equation}{\value{theorem}}
  }

\newcommand\enumcall[1]{
  \stepcounter{equation}
  \begin{#1}
  }

\newcommand\exitenumcall[1]{
  \end{#1}
  }

\newcommand\bdry{\partial}
\DeclareMathOperator{\cl}{cl}
\DeclareMathOperator{\interior}{int}
\DeclareMathOperator{\LH}{LH}  

\DeclareMathOperator{\ord}{ord}
\newcommand\card{\text{card}}  

\newcommand\union{\mspace{0.5mu}\mathbin{\cup}\mspace{0.5mu}}  
\newcommand\with{\mspace{2mu} {\circ} \mspace{2mu}} 

\newcommand\bbz{\mathbb{Z}}
\newcommand\bbr{\mathbb{R}}
\newcommand\bbc{\mathbb{C}}
\newcommand\bbh{\mathbb{H}}

\newcommand\eitheta{e^{\imath \theta}}

\newcommand\zbar{\bar z}

\newcommand\dx{d x}
\newcommand\dy{d y}

\newcommand\dmu{d \mu}
\newcommand\dnu{d \nu}
\newcommand\dtheta{d \theta}

\newcommand\poly[2]{P[#1,#2]}

\newcommand\abs[1]{\ensuremath{\vert #1 \vert}}
\newcommand\norm[1]{\ensuremath{\Vert #1 \Vert}}
\newcommand\inprod[2]{\ensuremath{\langle #1, #2 \rangle}}

\newcommand\plusperp{\ensuremath{\oplus_\perp}}

\newcommand\orthcomp{\ensuremath{\ominus_\perp}}

\begin{document}

\date{}

\title[An Abstract Beurlings Theorem II]{An Abstract Beurlings Theorem for Several Complex Variables II}
\author{Charles W. Neville}

\address{1 Chatfield Drive, Apartment 136\\
West Hartford, CT 06110\\
USA}

\email{chip.neville at gmail.com}


\thanks{Dedicated to the memory of the late Lee Rubel.  Lee, this one's for you.}

\begin{abstract}
In a previous paper, we presented an Abstract Beurling's Theorem for valuation Hilbert modules over valuation algebras. In this paper, we shall apply this theorem to obtain complete descriptions of the closed invariant subspaces of a number of Hilbert spaces of analytic functions in several complex variables.
\end{abstract}

\maketitle

\setcounter{section}{0}

\section{Introduction and notation.}\label{sec1}
In a previous paper, we stated and proved an Abstract Beurling's Theorem for valuation Hilbert modules over valuation algebras \cite[Theorem 8.1]{N3}. As an example of the power of  this theorem, we also used it to obtain a complete description of the closed invariant subspaces of $H^2$ of the polydisk.

Henceforth, we shall refer to our previous paper, \cite{N3}, as paper I, and we shall refer to this paper as paper II.

In this paper, paper II, we shall apply theorem 8.1 and its corollaries from paper I to obtain complete descriptions of the closed invariant subspaces of several Hilbert spaces of analytic functions in several complex variables. To do this, we first must review our standing notation from paper I. We shall quote almost verbatim, but without quotation marks:

 \
 
$\bbc$ will be the complex numbers and $\bbc^n$ complex n-space.  $\bbr$ will be the real numbers and $\bbr_+$ the positive reals, $\{r \in \bbr : r \geq 0\}$.  $\bbz$ will be the integers and $\bbz_+$ the positive integers, $\{m \in \bbz : m \geq 0\}$.  $E$ will be a complex Hilbert space.  $\poly{\bbc}{z_1, \dots, z_n}$ will be the algebra of polynomials with complex coefficients in the $n$ complex variables $z_1, \dots, z_n$, and $\poly{E}{z_1, \dots, z_n}$ will be the complex vector space (and $\poly{\bbc}{z_1, \dots, z_n}$ module) of polynomials with $E$-valued coefficients in the $n$ complex variables $z_1, \dots, z_n$.  We shall shamelessly identify polynomials with functions, so $\poly{E}{z_1, \dots, z_n}$ can be thought of as the space of $E$-valued polynomials in $z_1, \dots, z_n$.

We shall denote the Hilbert space norm on $E$ by $\norm{\cdot}_E$ and the inner product by $\inprod{\cdot}{\cdot}_E$.

$R$ will be a complex algebra (sometimes without a unit element) and $\bbh$ a left Hilbert module over $R$.  That is $\bbh$ will be both a Hilbert space and a left $R$ module, and the multiplication maps $m_r \colon h \mapsto rh$ will be continuous for each fixed $r \in R$.  We shall denote the Hilbert space norm on $\bbh$ by $\norm{\cdot}$ and the inner product by $\inprod{\cdot}{\cdot}$.

In the frequently occurring case where $\bbh$ is a Hilbert module of $E$-valued analytic functions, the norm and inner product on $E$ will be subscripted as indicated above, but they will be easy to distinguish from the corresponding quantities on $\bbh$, which will be unsubscripted.

If $X$ is a subset of a larger topological space $Y$, we shall denote the closure of $X$ by $\cl X = \cl_Y X$, the interior of $X$ by $\interior X = \interior_Y X$, and the topological boundary of $X$ by $\bdry X = \bdry_Y X$.  Finally, $\card(X)$ will denote the cardinality of the set $X$.

If $p \in \bbc^n$ and $0 < r < \infty$, $D^n(p, r)$ will denote the $n$-dimensional polydisk of radius $r$ centered at $p$, $D^n(p, r) = \{z \in \bbc^n \colon \abs{z_j - p_j} < r, j = 1 \dots n\}$. Similarly if $r = r_1 \dots r_n$ with $0 <  r_1 \infty, \dots 0 < r_n < \infty$ is a polyradius, $D^n(p, r)$ will denote the $n$-dimensional polydisk centered at $p$ of polyradius $r$.
 
\section{Valuation algebras and valuation {H}ilbert modules.}\label{sec2}

We also must review the definitions and basic properties of valuation algebras and valuation Hilbert modules from paper I. Again, we quote directly from that paper:

\thmcall{definition}\label{def2.1}
Let $R$ be a (possibly non-commutative) algebra.  An \emph{algebra valuation} on $R$ is a function $\ord \colon \! R \mapsto \bbz_+ \union \{\infty\}$ such that for all $r$ and $s \in R$,
\enumcall{enumerate}
  \item $\ord(r) = 0$ if $r$ is a left or right unit of $R$. \label{item2.1.1}
  \item $\ord(r) = \infty$ if and only if $r = 0$. \label{item2.1.2}
  \item $\ord(r s) \geq \ord(r) + \ord(s)$.  \label{item2.1.3}
  \item $\ord(\lambda r) = \ord(r)$ for $\lambda \in \bbc, \lambda \neq 0$ \label{item2.1.4}
  \item $\ord(r + s) \geq \min(\ord(r), \ord(s))$  \label{item2.1.5}
\exitenumcall{enumerate}
\exitthmcall{definition}
Of course, condition (\ref{item2.1.1}) is satisfied vacuously if $R$ does not have a two sided identity element.

Here and throughout, we follow the usual conventions with respect to $\infty$: $m < \infty$ for all $m \in \bbz$, and for such $m$,  $m \cdot \infty = \infty \cdot m = \infty$ if $m \neq 0$.  Further, $\infty \cdot 0 = 0 \cdot \infty = \infty$ and $\infty \cdot \infty = \infty$.

\thmcall{definition}\label{def2.2}
A \emph{valuation algebra} is an ordered pair $(R, \ord)$, where $R$ is a (possibly non-commutative) algebra and $\ord$ is an algebra valuation on $R$.
\exitthmcall{definition}

\thmcall{definition}\label{def2.3}
Let $(R, \ord_R)$ be a valuation algebra, and let $\bbh$ be a Hilbert space which is a left $R$ Hilbert module.  A \emph{Hilbert module valuation} on $\bbh$ (with respect to $(R, \ord_R)$) is a function $\ord \colon \bbh \mapsto \bbz_+ \union \{\infty\}$ such that for all $h$, $h_1$ and $h_2 \in \bbh$, and for all $r \in R$,
\enumcall{enumerate}
  \item $\ord(h) = \infty$ if and only if $h = 0$. \label{item2.3.1}
  \item $\ord(r h) \geq \ord_R(r) + \ord(h)$  \label{item2.3.2}
  \item $\ord(\lambda h) = \ord(h)$ for $\lambda \in \bbc, \lambda \neq 0$. \label{item2.3.3}
  \item $\ord(h_1 + h_2) \geq \min(\ord(h_1), \ord(h_2))$. \label{item2.3.4}
  \item The $\ord$ function is upper semi-continuous on $\bbh$. \label{item2.3.5}
\exitenumcall{enumerate}
\exitthmcall{definition}

Note that each set $\bbh_m= \{ h \in \bbh : \ord(h) \geq m \}$ is a subspace by properties (\ref{item2.3.3}) and (\ref{item2.3.4}).  The subspaces $\bbh_m $ are closed for $m = 0, 1, 2, \, \dots$ because $\ord$ is upper semi-continuous. The subspace $\bbh_0$ is automatically a closed subspace because $\bbh_0 = \bbh$. Fiinaly, as part of the definition of an $R$ Hilbert module, the multiplication map $m_r \colon h \mapsto r h$ is continuous for each $r \in R$.

\thmcall{definition}\label{def2.4}
Let $(R, \ord_R)$ be a valuation algebra.  A \emph{valuation Hilbert module} over $(R, \ord_R)$ is an ordered pair $(\bbh, \ord)$, where $\bbh$ is a left Hilbert module over $R$ and $\ord$ is a Hilbert module valuation on $\bbh$ with respect ot $(R, \ord_R)$.
\exitthmcall{definition}

Where it will cause no confusion, we shall denote both the algebra valuation on $R$ and the Hilbert module valuation on $\bbh$ by the same symbol, $\ord$.  We shall also make the gloss, whenever convenient, of denoting the valuation algebra $(R, \ord)$ by $R$ alone, and the valuation Hilbert module $(\bbh, \ord)$ by $\bbh$ alone.

\

To state our abstract Beurling's theorem, we need to review still more notation, conventions, and definitions from paper I:

Throughout, unless otherwise stated, $V$ shall be a closed subspace of $\bbh$.  To avoid trivial issues, unless otherwise stated, $V \neq \{ 0 \}$.

A valuation on a valuation algebra $R$ or a valuation Hilbert module $\bbh$ is \emph{strict} if equality holds in properties \ref{item2.1.3} or \ref{item2.3.2}. If the algebra $R$ or the left Hilbert module over $R$ admits a strict valuation, then the algebra must be quite well behaved. $R$ then has no left or right divisors of zero, and if $R$ is commutative, then it is an integral domain.

In definition 5.1 of paper I, the \emph{valuation subspace series} is defined to be the decreasing sequence of closed subspaces
  $V = V_0 \supseteq V_1 \supseteq V_2 \supseteq \cdots$, where $V_k = \{ h \in V : \ord(h) \geq k \}$. Similarly, $R_k = \{ r \in R \colon \ord{r} \geq k$.

In paper I, equation 3.12 ff., a closed subspace $W$ of $\bbh$ is defined to be \emph{near homogeneous} if the $\ord$ function is constant on all non-zero elements of $W$.

In paper I, definition 5.3, the \emph{near homogeneous decomposition} of $V$ is defined to be the orthogonal direct sum
\begin{equation*}
  V = W_0 \plusperp W_1 \plusperp W_2 \plusperp \, \cdots ,
\end{equation*}
where
\begin{equation*}
  W_0 = V_0 \orthcomp V_1, \; W_1 = V_1 \orthcomp V_2, \; W_2 = V_2 \orthcomp V_3, \, \dots \text{ .}
\end{equation*}
The orthogonal subspaces $W_0$, $W_1$, $W_2, \, \dots$ are called the \emph{near homogeneous components} of $V$.

In paper I, definition 6.3, the near homogeneous decomposition of $V$ is said to be \emph{near inner} if for each $k = 0, 1, 2 \dots$, $m = 0, 1, 2 \dots$, $h \in W_k$, $g \in V$, $r \in R_1$, and $\ord(r h + g) > m$, we have $ r h + g \perp W_m$.

In paper I, definition 7.1, the near homogeneous decomposition of $V$ is defined to have the \emph{full projection property} if for each $k$ and $m = 0, 1, 2, \dots$,
\begin{equation*}
  P_m^\bbh (r h + g) \in P_m^\bbh (W_m) ,
\end{equation*}
whenever $r \in R_1$, $h \in W_k$, $g \in V$, and $\ord(rh + g) \geq  m$.

\

Now we may restate our abstract Beurling's theorem and its most important corollary \cite[theorem 8.1 f.f.]{N3} from Paper I. As shown in Paper I, this theorem and its corollaries completely describe the closed invariant subspaces of valuation Hilbert modules over valuation algebras.

\thmcall{theorem}\label{thm2.5}
Let $V$ be a closed subspace of a valuation Hilbert module over a valuation algebra $R$.  Then $V$ is $R_1$ invariant if and only if the near homogeneous decomposition of $V$ is near inner and has the full projection property.
\exitthmcall{theorem}

\thmcall{corollary}\label{cor2.6}
Let $R$ be a valuation algebra satisfying either of the following conditions:
\enumcall{enumerate}
  \item $R = R_1$, or
  \item $R$ is unital
\exitenumcall{enumerate}
Then a closed subspace $V$ of a valuation Hilbert module $\bbh$ over $R$ is an $R$ submodule of $\bbh$ if and only if the near homogeneous decomposition of $V$ is near inner and has the full projection property.
\exitthmcall{corollary}

\section{Analytic {H}ilbert modules.}\label{sec3}

In this section, we shall introduce the general class of concrete valuation algebras and valuation Hilbert modules to which our theorem \ref{thm2.5} and its corollary apply.

\thmcall{definition}\label{def3.1}
Let $\Omega$ be a connected paracompact analytic manifold of dimension $n \geq 1$, with a  distinguished basepoint $p \in \Omega$. Let $E$ be a complex Hilbert space. 
 
Let $R$ be a complex algebra of $\bbc$ valued analytic functions on $\Omega$.  Define the order of $f \in R$, $\ord(f)$, to be the order of the zero of $f$ at $p$.  Let $\bbh$ be a complex Hilbert space of $E$ valued analytic functions on $\Omega$. Define the the order of $h \in \bbh$, $\ord(h)$, to be the order of the zero of $h$ at $p$.

$\bbh$ is an analytic Hilbert module on $\Omega$ over the analytic algebra $R$ if \,$\bbh$ is a left Hilbert module over $R$ under pointwise operations, and if each sequence $h_n, n = 1, 2, 3, \dots$ in $\bbh$ converging to $h$ in $\bbh$ norm also converges to $h$ uniformly on compact subsets of $\Omega$. In other words, if the Hilbert space topology on $\bbh$ is finer than the topology of compact convergence on $\bbh$.
\exitthmcall{definition}

Note that our analytic Hilbert modules are different from Guo's analytic Hilbert modules \cite{Guo1}.  Note also that the $\ord$ function is invariant across all local coordinate systems in $\Omega$.

\thmcall{theorem}\label{thm3.2}
Let $\bbh$ be an analytic Hilbert module on $\Omega$ over the analytic algebra $R$. Then $\bbh$ is a valuation Hilbert module over the valuation algebra $R$. Thus theorem \ref{thm2.5} and its corollary completely describe the closed $R_1$ and $R$ invariant subspaces of $\bbh$.
\exitthmcall{theorem}

\begin{proof}
All we have to show is that the $\ord$ function on $\bbh$ is upper semi-continuous at $p$.  To this end, let $V$ be an open coordinate neighborhood of the basepoint $p$ in $\Omega$, and let $\phi \colon V \rightarrow \phi(V)$ be a local coordinate mapping.  Note for clarity that $\phi(V) \subseteq \bbc^n$. Without loss of generality, $\phi(p) = 0$.

To simplify notation, let $\psi = \phi^{-1}$.  Let $0 < r < 1/2$. and let $D^n(0, 2r)$ be an $n$-dimensional polydisk centered at $0$ such that $\psi(\cl D^n(0, 2r)) \subset V$. Note that $\psi(0) = p$.

Now let $h_n, n = 1, 2, 3, \dots$be a norm convergent sequence in $\bbh$ converging to $h \in \bbh$.  By definition \ref{def3.1}, the sequence $h_n$ also converges uniformly on compact subsets of $\Omega$ to $h$.  Clearly, $\psi(\cl D^n(0, 2r))$ is such a subset, so the sequence $h_n$ converges uniformly on $\psi(\cl D^n(0, 2r))$ to $h$. To simplify notation, let $f_n = h_n \with \psi$, and $f = h \with \psi$.  Then the sequence $f_n$ converges uniformly on $\cl D^n(0, 2r)$ to $f$.

From definition \ref{def3.1}, we see that $f$ is an ordinary analytic function, and $\ord(h) = \ord(f)$ is the common value of the absolute value of the multi-indices of the first non-vanishing terms of the power series expansion of $f$ at $0$,  or $\infty$ if $f$ is identically zero. In other words, if $a_{m_1} \dots a_{m_n} z^{m_1} \cdots z^{m_n}$ is one of these terms, then $\ord(h) = \ord(f) = \abs{m_1} + \cdots + \abs{m_n}$.

On the polydisk $D^n(0, r)$, 
\begin{equation*}
   r^m a_m = \int_0^{2\pi} f(r e^{\imath \theta}) e^{-\imath m \theta} \dtheta/(2 \pi^n)
\end{equation*}
in multi-index notation, with a corresponding formula for $r^m a_{j, m}$ with $f$ replaced by $f_j$.  Here, the integral is repeated $n$ times, $\theta = \theta_1 \cdots \theta_n$, and $\dtheta = \dtheta_1 \cdots \dtheta_n $.

Since the sequence $f_j$ converges uniformly to $f$ on the distinguished boundary of the polydisk $D^n(0, r)$, the power series coefficients $a_{j, m}$ of $f_j$ converge to $a_m$ as $j \rightarrow \infty$.  Now $\ord(f_j) = \abs{m_j}$ and $\ord(f) = \abs{m}$, so by the usual argument, $\liminf{ \ord(f_j) } \geq \ord(f)$. Hence the $\ord$ function is upper semi-continuous at $0$ on $\{ h \with \psi \colon h \in \bbh \}$, and thus the $\ord$ function is upper semi-continuous at $p$ on $\bbh$.
\end{proof}

\section{$H^2$ spaces.}\label{sec4}
Let $\Omega$ be a domain in $\bbc^n$ rather than a general connected analytic manifold, and let $R = \poly{\bbc^n}{z_1, \dots, z_n}$ or $H^\infty(\Omega)$.

\

In paper I, we discussed $H^2(E, D^n(0, 1)))$ of the unit polydisk centered at $0$. We showed that our theorem \ref{thm2.5} and its corollary (i{.}e{.}~\cite[theorem 8.1]{N3} and its corollaries in paper I) completely characterize the closed invariant subspaces of this Hilbert module.  By our discussion in the last section, section \ref{sec3}, this result carries over to $H^2(D^n(p, r))$, where  $0 < r < \infty$ and $p \in \Omega$.  However, there are many more domains in  $\bbc^n$ in addition to polydisks on which we can define $H^2$ spaces.

\

Let $0 < r < \infty$, let $B^n(a, r)$ be the open ball of radius $r$ in $\bbc^n$ centered at $a$, and let $M$ be the group of analytic automorphisms of $B^n(a, r)$. Rudin essentially defined $H^p$ of the open unit ball in $\bbc^n$ to be the space of all analytic functions $h$ on $B^n(0, 1)$ such that the $M$-subharmonic function $z \rightarrow \abs{h(z}^p$ has an $M$-harmonic majorant. \cite[definition 5.6.1 pp 83--84 ]{Ru2}  Here, $M$ is the group of conformal automorphisms of the unit ball, and an $M$-harmonic function is a function annihilated by the $M$ invariant Laplacian \cite[sections 4.1--4.4]{Ru2}. We shall adopt Rudin's definition essentially without change.

\thmcall{definition}\label{def4.1}
   Let $0 < r < \infty$ and let $p \in \bbc^n$.  Define $H^2(E, B^n(p, r))$ to be the Hilbert space of all analytic functions $h$ on $B^n(p, r)$ such that the $M$-subharmonic function $z \rightarrow \norm{h(z)}_E^2$ admits an $M$- harmonic majorant.
   
If $h$ is an analytic function in $H^2(E, B(p, r))$, let $u$ be the above least $M$-harmonic majorant. We define the $H^2(E, B^n(p, r))$ norm of $h$ to be $\norm{h}^2 = u(p)$.
\exitthmcall{definition}

For those more comfortable with harmonic functions, Rudin observes that $M$-harmonic functions and $M$-subharmonic functions may be replaced by ordinary harmonic functions and subharmonic functions in the definition of $H^2(E, D^n(0, 1))$ without changing anything essential, except that then it is not obvious that \break
$H^2(E, D^n(0, 1))$ is $M$ invariant.

\thmcall{theorem}\label{thm4.2}
$H^2(E, B^n(p, r))$ is an analytic Hilbert module over $R$. Therefore, by theorem \ref{thm3.2}, our abstract Beurling's theorem, theorem \ref{thm2.5}, and its corollary completely describe the closed invariant subspaces of $H^2(E, B^n(p, r))$.
\exitthmcall{theorem}

\thmcall{proof}
This will follow quickly from the Harnack inequalities for $M$-harmonic functions. The Harnack inequalities themselves follow immediately from the formula for the invariant Poisson kernel \cite[definition 3.3.1]{Ru2}.

Let $0 < r < \infty$. Let $\Omega = B^n(a, r)$ and let $a \in \Omega$. To simplify notation, let $\bbh = H^2(E, \Omega)$.  Recall that if $h \in \bbh$, $\norm{h}$ is the norm of the function $h$ in $\bbh$, and if $z \in \Omega$, then $\norm{h(z)}_E$ is the norm of the element $h(z)$ in $E$.

Let $h_j, j = 1, 2, 3, \dots$ be a norm convergent sequence in $\bbh$ converging to $h \in \bbh$. Let $f_j = h - h_j$. Then  $\norm{f_j} \rightarrow 0$ as $j \rightarrow \infty$.

Let $v_j$ be the least $M$-harmonic majorant of the $M$-subharmonic function $z \rightarrow \norm{f_j(z)}^2_E$.  Let $K$ be a compact subset of $\Omega$.  By the Harnack inequalities, there is a constant $0 < C < \infty$ such that $v_j(w) < C v_j(a)$ for all $w \in K$.  But $\norm{f_j(w)}_E^2 \leq v_j(w) \leq C v_j(a) = C \norm{f_j}^2$, so $f_j$ converges uniformly on the compact set $K$ to $0$. Thus $h_j$ converges uniformly to $h$ on compact subsets of $\Omega$.
\exitthmcall{proof}

Both $D^n(p,r) $ and $B^n(p,r)$ are examples of \emph{bounded symmetric domains}.  And, in fact, if $\Omega$ is a bounded symmetric domain with basepoint $p$,, we can define $H^2(E, \Omega)$ in the same way as for $H^2(E, D^n(p,r), p)$ and $H^2(E, B^n(p,r))$.  However, we have to make a distinction already present in our treatments of $H^2(E, B^n(p,r))$ and $H^2(E, D^n(p,r))$:

\

If $\Omega$ is \emph{irreducible} (not the product of bounded symmetric domains of lower dimension), we define $H^2(E, \Omega)$ to be the Hilbert space of all $E$-valued analytic functions $h$ on $\Omega$ such that the $M$-subharmonic function $ s(z) = \norm{h(z)}_E^2$ has an $M$-harmonic majorant, exactly as in the case of the ball $B^n(a,r)$.  Here, as in the case of the ball, $M$ is the group of automorphisms of $\Omega$.  And, exactly as in the case of the ball, if $u_h$ is the least harmonic majorant of $s_h$, we norm $H^2(E, \Omega)$ by $\norm{h}^2 = u_h(p)$.  Though the norm depends on the choice of the basepoint $p$, by the Harnack inequalities, different choices of $p$ yield equivalent norms.

\

If, on the other hand, $\Omega$ is \emph{reducible} (the product of $k$ irreducible bounded symmetric domains $\Omega_1 \times \cdots \times \Omega_k$ of lower dimension), we define $H^2(E, \Omega)$ to be the Hilbert space of all $E$-valued analytic functions $h$ on $\Omega$ such that each $j$-$M_j$-subharmonic function $s_{h_j}(z) = \norm{h_j(z)}_E^2$ has a $j$-$M_j$-harmonic majorant, exactly as in the case of the polydisk $D^n(a, r)$.

In the above, $M_j$ is the group of automorphisms of $\Omega_j$, and  $h_j$  is the slice function $z_j \rightarrow h(a_1, \dots, z_j, \dots a_k)$ on $\Omega_j$. The function $h$ is $k$-$M$-harmonic if it is continuous at each $a = (a_1, \dots, a_n)$, and each slice function $h_j$ is $M_j$-harmonic near $a_j$.  Here, as one might expect, $a_1, \dots, a_k \in \Omega_k$ and $z_j \in \Omega_j$.  If $u_h$ is the least $k$-$M$-harmonic majorant of $s_h$, we norm $H^2(E, \Omega)$ by $\norm{h}^2 = u_h(p)$, exactly as in the case of the polydisk.  Again, though the norm depends on the choice of the basepoint $p$, by the Harnack inequalities, different choices of $p$ yield equivalent norms.  In summary, we have

\thmcall{theorem}\label{thm4.3}
Let $\Omega$ be a bounded symmetric domain with basepoint $p \in \Omega$ Then $H^2(E, \Omega)$ is an analytic Hilbert module over the analytic algebra $R$. Therefore, by theorem \ref{thm3.2}, our abstract Beurling's theorem and its corollary completely describe the closed invariant subspaces of $H^2(E, \Omega)$.
\exitthmcall{theorem}

\thmcall{proof}
See the above.
\exitthmcall{proof}

\thmcall{exercise}\label{exercise4.4}
If you are unfamiliar with bounded symmetric domains, look up the definition and show that $B^n(a, r)$ is an irreducible bounded symmetric domain, while $D^n(a, r)$ is a reducible bounded symmetric domain.  To see proofs that $D^n(a,r)$ and $B^n(a, r)$ are bounded symmetric domains, but not using the that terminology, read Rudin, specifically the beginning chapters of \cite{Ru1} and \cite{Ru2}.
\exitthmcall{exercise}

An extremely general example of an analytic Hilbert module is the case of $\bbh = H^2(E, \Omega, p)$, where $\Omega$ is a connected complex analytic manifold, the basepoint $p \in \Omega$, , and $R = H^\infty(\bbc, \Omega)$.

Unfortunately, there are a number of difficulties in this case, beginning with the Laplacian on $\Omega$ (the Laplace Beltrami operator), which in general has variable coefficients.  Thus we have to appeal to the modern theory of PDEs to establish the potential theory we need.  Second, we need a Hermitian structure on $\Omega$ before we can even define the Laplace Beltrami operator.  Third, we need to require that $\Omega$ is a K\"{a}hler manifold (or compact, the trivial case) in order to guarantee that analytic functions are harmonic.  So we shall defer the discussion of Hardy class $H^2$ on complex analytic manifolds to another paper.

\

There is, however, one circumstance where the difficulties mentioned above are already worked out, and that is the the case where the the (complex) dimension of $\Omega$ is $1$, in other words where $\Omega$ is a Riemann surface with basepoint $a$..  In fact, this is the case where Walter Rudin first introduced the harmonic majorant definition of $H^2$, and of $H^p$ for $1 \leq p < \infty$ \cite{Ru3}.  While the ordinary Laplacian is not invariant under coordinate changes, ordinary functions harmonic in local coordinates on $\Omega$ are.  Furthermore, analytic functions are harmonic in local coordinates.  The passage to $H^2(E, \Omega)$ is immediate by taking component functions.  Thus, if we define $H^2(E, \Omega)$ as the space of $E$-valued analytic functions $h$ on $\Omega$ such that the subharmonic function $s_h(z) = \norm{h(z)}_E^2$ has a harmonic majorant, and if we norm $H^2(E, \Omega)$ by $\norm{f}^2 = u_h(a)$, where $u_h$ is the least harmonic majorant of $s_h$, we have, just as in the cases already treated of the ball and bounded symmetric domains,

\thmcall{theorem}\label{thm4.5}
Let $\Omega$ be a Riemann surface with basepoint  $a \in \Omega$, $\bbh = H^2(E, \Omega)$, and $R = H^\infty(\bbc, \Omega)$,  Then $\bbh$ is an analytic Hilbert module over the analytic algebra $R$. Therefore, by theorem \ref{thm3.2}, our abstract Beurling's theorem and its corollary completely describe the closed invariant subspaces of $\bbh$.
\exitthmcall{theorem}

What is astounding about this example is its generality.  The author in \cite{N1} and \cite{N2} presented a different description of the closed submodules of $H^2(\bbc, \Omega)$ in terms of multiple valued inner functions $I$ with a single valued absolute value $\abs{I}$.  (In fact, the author described the closed submodules of $H^p(\bbc, \Omega)$ for $1 \leq p < \infty$ and the $\beta$ closed ideals of $H^\infty(\bbc, \Omega)$.)  Later, Hasumi, \cite{Has1} and \cite{Has2}, extended this description to $L^p(\bbc, \dmu)$ for $1 \leq p < \infty$, where $\dmu$ is harmonic measure on the Martin boundary of $\Omega$. Hasumi used the language of analytic vector bundles rather than multiple valued analytic functions.

In the above, the Riemann surface, while it could be infinitely connected, had to satisfy strict abstract or geometrical conditions.  The author in his thesis, using unpublished work of Kennel, was able to give an example of an infinitely connected \emph{plane} domain where the theorem fails for $H^\infty(\bbc, \Omega)$ under the $\beta$ topology.  (See \cite{N2} for a summary.)  By contrast, theorem \ref{thm4.5} applies to Hilbert space valued functions and to \emph{all} Riemann surfaces.

It is rare that results in \emph{several} complex variables, such as theorem \ref{thm3.2} and theorem \ref{thm4.5}, combined with our abstract Beurling's theorem, give improved results in \emph{one} complex variable, but this is one of those rare cases.

\section{$A^2$ spaces.}\label{sec5}
Let $\Omega$ be a domain in $\bbc^n$, and let $\dmu$ be $2 n$-dimensional volume measure on $\bbc^n$. The Bergman space $A^2(E, \Omega, \dmu)$ is the space of all $E$-valued analytic functions $h$ on $\Omega$ such that

\begin{equation*}
    \int_\Omega \norm{h(z)}_E^2 \dmu (z) < \infty .
\end{equation*}

A celebrated result of Aleman, Richter, and Sundberg from the 1990s established a Beurling type theorem characterizing the closed $\poly{\bbc}{z}$ submodules of $A^2(\bbc, D^1(0, 1), \dmu)$ \cite{ARS1}, where $\dmu$ is real $2$-dimensional area measure.  Their proof was quite difficult, and famously involved biharmonic functions.  Since then, considerably simpler proofs have been given by Shimorin, and later Izuchi, and Sun and Zheng (\cite{Sh1}, \cite{Sh2}, \cite{Iz1}, \cite{SZ1}).  Shimorin has extended the ARS theorem to weighted $A^2$ spaces, $A^2(\bbc, D^1(0, 1), k\dmu)$ of the unit disk, albeit with restrictions on the weights $k$ such as circular symmetry and a subharmonicity condition (\cite{Sh1}, \cite{Sh2}).  And recently, Ball and Bolotnikov gave an elegant and surprising characterization of the closed $\poly{\bbc}{z}$ submodules of weighted Hilbert space valued $A^2$ spaces of the unit disk (with the particular weight functions $(1 - \abs{z}^2)^{k - 2}$ for $k \geq 2$) \cite{BaBo1}.

\

In this section, we shall generalize the ARS and Shimorin theorems to several complex variables and to connected paracompact analytic manifolds.

\

Let $\Omega$ be a connected, paracompact analytic manifold with basepoint $p$.  Let $\dmu$ be a real $2 n$-dimensional volume measure on $\Omega$. Let $k$ be an admissible weight function on $\Omega$.  (We shall discuss the existence of the measure $\dmu$ and give the definition of an admissible weight function momentarily.)  Let $R = H^\infty(\Omega)$,  the algebra of bounded complex valued analytic functions on $\Omega$, or, in the case where $\Omega \subseteq \bbc^n$, let $R = \poly{\bbc}{z_1, \dots, z_n}$.

\thmcall{definition}\label{def5.1}
The $E$ valued weighted Bergman space on $\Omega$, $\bbh = A^2(E, \Omega, k\dmu)$, is the space of all $E$ valued analytic functions $h$ on $\Omega$ such that the weighted $\bbh$ norm
  \begin{equation*}
     \norm{h}^2 = \int_\Omega \norm{h(z)}_E^2 k(z) \dmu (z) < \infty .
  \end{equation*}
\exitthmcall{definition}

Define the order functions $\ord(f)$ and $\ord(h)$ in the usual way, as the order of the zero of $f \in R$ or $h \in \bbh$ at $p$, respectively.

\thmcall{theorem}\label{thm5.2}
    $\bbh = A^2(E, \Omega, k \dmu)$ is an analytic Hilbert module over the analytic algebra $R$.  Thus by our theorem on analytic Hilbert modules, theorem \ref{thm3.2}, our abstract Beurling's theorem, \ref{thm2.5} and its corollary, completely characterize the closed $R$ invariant subspaces of $\bbh$.
\exitthmcall{theorem}

In the special case of unweighted Bergman spaces defined on domains in $\bbc^n$, theorem \ref{thm5.2} follows immediately from the following lemma:

\thmcall{lemma}\label{lemma5.3}
  Let $V$ be a domain in $\bbc^n$.  Let $\bbh = A^2(\bbc, V, \dmu))$, and let $K$ be a compact subset of $V$. Then there is a constant $0 < C_K < \infty$ depending on $V$ and $K$, such that
   \begin{equation*}
     \sup_{z \in K} \abs{h(z)} \leq C_K \norm{h} \text{   all  }  h \in A^2 (\bbc, V, \dmu) \text{  (\cite{Kr1}, lemma 1.4.1)}.
   \end{equation*}
Thus a sequence $h_k, k \in \bbz_+$ converging in norm to $h$ in $\bbh$ also converges to $h$ in the topology of compact convergence on compact subsets of $V$.  Hence $\bbh$ is  an analytic Hilbert module over the analytic algebra $R$.
\exitthmcall{lemma}  

The passage to $\bbh = A^2(E, V)$ is immediate.

\

For the special case of the unit ball in $\bbc^n$, Rudin calculates the Bergman kernel (Rudin calls it the Cauchy kernel) for  $A^2$. (\cite{Ru2}, definition 3.2.)  Although Rudin does not do this, his calculation allows for a proof of lemma \ref{lemma5.3} by estimating the Bergman kernel.

\

To return to volume measures and admissible weight functions, recall that a real $2n$-dimensional volume measure on a connected paracompact complex analytic manifold $\Omega$ is a measure which, in local coordinates $z = z_1, \dots, z_n$ on a typical open coordinate neighborhood $V$, is $\dmu(z) = c(z) \dx_1 \dy_1 \cdots \dx_n \dy_n $.  Here, $c$ is Lebesgue measurable (with respect to the measure $\dx_1 \dy_1 \cdots \dx_n \dy_n$ on $V$), and bounded {a.}{e.}~on $V$. Further, for some open coordinate neighborhood $V$ of the basepoint $p$, $c >$ some $\delta > 0$ on $V$ {a.}{e.}   A suitable volume measure always exists on $\Omega$, and in fact there are many of them.  For example, if $G$ is a Riemannian metric on $\Omega$, the Riemannian volume element $\dmu = \sqrt{ \text{det}(G)} \, \dx_1 \dy_1  \cdots  \dx_n \dy_n$ is such a measure \cite{Gr1}.

\thmcall{definition}\label{def5.4}
A real valued function $k$ is an admissible weight function on $\Omega$ if it is Lebesgue measurable and locally bounded away from $0$ a.e.\@ on $\Omega$.  To be more specific, if for each $a \in \Omega$,  there exists a connected open coordinate neighborhood $V$ of $a$, such that $k$ is bounded away from $0$ a.e.\@ on $V$,  that is if there exists a constant $\epsilon_V >  0$ such that $k > \epsilon_V$ a.e.\@ on $V$.
\exitthmcall{definition}

Note that our admissible weight functions are considerably more general than Shimorin's.

\begin{proof}[Proof of Theorem \ref{thm5.2}]
Let $V$ be an open coordinate neighborhood as in definition \ref{def5.4}, and $\phi \colon V \rightarrow U$ be an analytic chart. Note for clarity that $U$ is a domain in $\bbc^n$. Without loss of generality, $\phi(p) = 0$.  To simplify notation, let $\psi = \phi^{-1} \colon U \rightarrow V$.

In the notation of definition \ref{def5.4}, let $0 < \epsilon \leq k(z)$ a.e.\@ on $V$, and let $K$ be a compact subset of $V$. Let $C_{\phi(K)}$ be as in lemma \ref{lemma5.3}, with $K$ replaced by the compact subset $\phi(K)$ in $U$. Let $C_K = C_{\phi(K)}$.  The theorem now follows from the chain of inequalities,
\thmcall{equation*}
  \begin{split}
\epsilon \sup_{z \in K} \norm{h(z)_E^2} &\leq \epsilon \,C_K\int_U \norm{h \with \psi(z)}_E^2 \dnu(z) = \epsilon \,C_{K} \int_V \norm{h(z)}_E^2 k(z) \dmu(z) \\
          &\leq \epsilon \,C_K \int_\Omega \norm{h(z)}_E^2 k(z) \dmu(z) = \epsilon \,C_K \norm{h}_E^2 .
   \end{split}
\exitthmcall{equation*}

Here, the left hand inequality follows from lemma \ref{lemma5.3}, and the passage from the integral over $V$ to the integral over $\Omega$ follows from the monotoneity of the Lebesgue integral for positive functions.

Thus the topology of compact convergence is finer than the Hilbert space topology on $\bbh$, and so $\bbh$ is an analytic Hilbert module over the analytic algebra $R$.
\end{proof}

As we observed at the end of section \ref{sec4}, it is rare for a result in several complex variables to generalize results on one complex variable, but this is one of those rare cases.

\section{Paradox lost.}\label{sec6}
The title refers, not to Milton's great poem, but to an apparent paradox which arises from the Aleman, Richter, and Sundberg (ARS) characterization of the closed invariant subspaces of the (unweighted) Bergman space of the unit disk in $1$ complex variable, that is $A^2( \bbc, D^1(0, 1), d\mu)$ in our notation.

\thmcall{theorem}[Aleman, Richter, and Sunberg \cite{ARS1}]\label{thm6.1}
 Let $V$ be a closed invariant subspace of $\bbh = A^2(\bbc, D^1(0, 1), \dmu)$, where $\dmu$ is real $2$-dimensional area measure. Then $V$ is the smallest closed invariant subspace generated by $M$\!, where $M = V \orthcomp z V$\!.
\exitthmcall{theorem}

The ARS theorem is different, and arguably better, than ours. The apparent paradox it raises involves the dimension of $M$. Let $V = \plusperp W_m$ be the near homogeneous decomposition of $V$\!. Recall that by proposition 5.17 of paper I, \cite{N3}, the projection operator
\begin{equation*}
    L_{ W_m }  \colon W_m \mapsto P_m^\bbh ( W_m ) \subseteq H_m
\end{equation*}
is invertible.  Now each subspace $H_m = \LH ( z^m )$ has dimension $1$. Thus each subspace
    \begin{equation*}
        W_m = V_m \orthcomp \! \! V_{m + }, \, \,  m =  1, 2, 3, \ldots
    \end{equation*}
has dimension $0$ (if it equals $\{ 0 \}$) or $1$. But, as Aleman, Richter, and Sunberg showed, $M$ can have any dimension from $1$ up to and including $\infty$.\\

The resolution of this paradox is contained in a general observation, which we shall state momentarily as a proposition. Recall from paper I, definition 5.1, that $V_k = \{ h \in V \colon \ord(h) \geq k \}$.  Note that if $k$ is the first index such that $W_k \neq \{ 0 \}$, then $V = V_{ k }$. Thus $M = V_{ k } \orthcomp z V_{k}$. Our general observation is,

\thmcall{proposition}\label{prop6.2}
    Let $ \bbh $ be an analytic Hilbert module over the analytic algebra $ P[\bbc, z] $. Let $ k $ be the first index index such that $ W_k \neq \{ 0 \} $. Suppose $ M $ has dimension $ \geq 2 $. Then $ z V_k $ is a proper subset of $ V_{ k+ 1} $, so $ M $ is a proper superset of $ W_{ k + 1 } $.
\exitthmcall{proposition}

\begin{proof}
    If $ h \in V_k $, then $ \ord( h ) \geq k $, so $ \ord ( zh ) \geq k+ 1 $. Thus $ z V_k \subseteq V_{k + 1} $, so $ M = V_k \orthcomp zV_k  \supseteq  V_k \orthcomp V_{ k + 1 } = W_k $. Now if $ M $ has dimension $ \geq 2 $, then $ M \supsetneq W_k $, so $ z V \subsetneq V_{ k + 1 } $.
\end{proof}

The difficulty with this proof is that it appears to be circular, in that it appears to use the paradox to resolve the paradox. Here is more direct proof for $A^2(\bbc, D^1(0, 1), \dmu)$, simply using the power series expansions of analytic functions in $A^2(\bbc, D^1(0, 1)), \dmu)$:

\begin{proof}
    As above, $ z V_k \subset V_{k + 1} $. Now suppose $ z V_k $ actually equals $ V_{ k + 1 } $. Let $ f $ and $ g $ be two linearly independent elements of $M$. Note that $ f $ and $ g \in V_k $. Let $ f ( z ) = a_0 + a_1 z + a_2 z^2 +\cdots $, and $ g ( z ) = b_0 + b_1 z + b_2 z^2 + \cdots $. There are two cases:
    
    If both of $ a_k $ and $ b_k \neq 0 $, then, by multiplying $ g $ by a constant if necessary, without loss of generality, $ a_k = b_k $. Let $ h = f - g $, and let $ h ( z ) = c_0 + c_1 z + c_2 z^2 + \cdots $. Then $ h \in M $ because $ M $ is a linear subspace, $ h \neq 0 $ because $ f $ and $ g $ are linearly independent, and $ c_k = a_k - b_k = 0 $, so $ h \in V_{ k + 1 } $.
    
    On the other hand, if $ a_k = b_k = 0 $, set $ h = f $, and note that, as above, $ h \in M $, $ h \neq 0 $, and $ c_k = 0 $, so $ h  \in V_{ k + 1 } $. Thus, in either case, there is a non-zero function $ h \in M $ such that $ h \in V_{ k + 1 } $. But since $ M = V_k \orthcomp V_{ k + 1 } $, $ h \perp h $, so $ h = 0 $. Thus $ h \neq 0 $, and $ h = 0 $, which is impossible.
\end{proof}

Hedenmalm, Richter, and Seip have given a more explicit construction of cases, involving the angle between subspaces, where $M = V \orthcomp zV$ has finite or infinite dimension. Their construction is valid, not only for $A^2(\bbc, D^1)$, but also for weighted Bergman spaces in the unit disk \cite{HRS1}.

Walter Rudin has given an example of a closed invariant subspace $V$ of $H^2(\bbc, D^2)$ (with $D^2 = D^2(0, 1) =$ the unit polydisk in $\bbc^2$ and $R = \poly{\bbc} {z_1, z_2}$) which is not finitely generated (\cite{Ru1}, Corollary to Theorem 4.4.1). While this not surprising, since $H_k$ is spanned by $z_1^k, z_1^{k-1}z_2, \cdots z_2^k$ and so has dimension $k$ in this case, it does show that being infinitely generated is as much the rule as the exception for analytic Hilbert modules in several complex variables.

\section{$A^2$ spaces of differentials.}\label{sec7}
The main point of this section is that it is possible for the elements of a valuation Hilbert module $\bbh$ to be something other than analytic functions. In this section, we shall show that they can be analytic differentials.

\

As before, let $\Omega$ be a connected paracompact analytic manifold of dimension $n \geq 1$, with a  distinguished basepoint $p \in \Omega$.  In multi-index notation, let $z = z_1, \dots, z_n$ and $w = w_1, \dots, w_n$ be local coordinate systems on $\Omega$.  Let $\alpha = f(z)\,dz_1 \wedge \cdots  \wedge dz_n$ be a complex valued analytic differential $n$-form on $\Omega$.  Then $\bar{\alpha} = \overline{f (z)}\,d \zbar_1 \wedge \cdots \wedge d \zbar_n$.  Note that $\bar{f}$ is anti-analytic, but so are the coordinate change coefficients from $\bar{z}$ to $\bar{w}$, $\partial \bar{w}_j / \partial \bar{z_j}$.  Thus $\bar{\alpha}$ is also invariant under local changes of variables. Consequently,
\begin{equation*}
    \begin{split}
          \alpha \wedge \bar{\alpha} &= f(z) \overline{ f(z) } dz_1 \wedge \cdots \wedge dz_n \wedge d\bar{z_1}  \wedge \dots \wedge d\bar{z}_n \\
       &= (-1)^{(n (n-1) / 2)} \,\abs{f(z)}^2 dz_1 \wedge d\bar{z}_1 \wedge \cdots \wedge dz_n . \\\wedge d\bar{z}_n        &= P_n \,\abs{f(z)}^2 dx_1 \wedge dy_1 \wedge \cdots \wedge dx_n \wedge dx_n ,
     \end{split}
\end{equation*}
where the constant $P_n = {(-2 \imath)}^n (-1)^{(n (n-1) / 2)}$. Consequently, $\alpha \wedge \bar{\alpha}$ is invariant under local coordinate changes. In summary,

\thmcall{lemma}\label{lemma7.1}
Let $\Omega$, $\alpha$, $f$, and $P_n$ be as above. Then $\alpha \wedge \bar{\alpha} / P_n$ is a positive real $2n$-dimensional, coordinate invariant volume measure on the under lying real manifold of $\Omega$, and
\begin{equation*}
  \alpha \wedge \bar{\alpha} / P_n = \abs{f(z)}^2 \,dx_1 dy_1 \wedge \dots \wedge dx_n \wedge dy_n .
\end{equation*}
\exitthmcall{lemma}

Let $ E$ be a complex Hilbert space.  To extend lemme \ref{lemma7.1} to $E$ valued analytic differentials, we need to define a conjugate linear involution $\alpha^{\!\star}$ for such differentials, analogous to $\bar{\alpha}$.  A conjugate linear involution on $E$ is an isometric mapping  $f \rightarrow f^\star$ from $E$ onto $E$, which is conjugate linear. In other words, $(a f + b g)^\star = \bar{a} f^\star + \bar{b} g^\star$, and which is an involution, that is $f^{\star \star} = f$.

Define $f \cdot g = \inprod{f}{g^\star}_{\!E}$. Note that $f \cdot g$ is a bilinear form on $E \times E$, and $f \cdot g^\star = \inprod{f}{g}_{\!E}$, the ordinary inner product of f and g on E.

Let $z = z_1, \dots z_n$ be a local coordinate system on $\Omega$. From now on $f$ shall be an $E$ an $E$ valued analytic function in this local coordinate system rather than an element of $E$.  Let $\alpha = f(z)\,dz_1 \wedge \cdots  \wedge dz_n$ be an $E$ valued analytic differential $n$-form on $\Omega$.  Define $\alpha^\star$ to be $\alpha^\star = f(z)^\star\,d\bar{z}_1 \wedge \cdots \wedge d\bar{z_n}$. Then $\alpha^\star$ is an $E$ valued anti-analytic differential $n$-form on $\Omega$. To see this, let $f(z) =\sum c_k z^k$ be the power series expansion of $f$ in, say, an open polydisk centered at $0$.. Then in multi-index notation,
\begin{equation*}
   \alpha^\star = f(z)^\star d\bar{z}_1 \wedge \cdots \wedge d\bar{z}_n = \sum \bar{c}_k \bar{z}^k d\bar{z_1} \wedge \cdots \wedge d\bar{z_n} , 
\end{equation*}
which is clearly anti-analytic.  As in the discussion preceding lemma \ref{lemma7.1}, so are coordinate change coefficients, so $\alpha^\star$ is invariant under local changes of variables.

We still need to define wedge products of $E$ valued $n$ forms.  In short, define
\begin{equation*}
  \begin{split}
   \alpha \wedge \alpha^\star &= f(z) \cdot f(z)^\star dz_1 \wedge \cdots \wedge dz_n \wedge d\bar{z_1}  \wedge \dots \wedge d\bar{z}_n \\
      &= \inprod{f(z)}{f(z)}^\star_E  dz_1 \wedge \dots \wedge dz_n \wedge d\bar{z}_1 \wedge \dots \wedge d\bar{z}_n \\
      &= (-1)^{(n (n-1) / 2)} \,\norm{f(z)}_E^2 dz_1 \wedge d\bar{z}_1 \wedge \cdots \wedge dz_n \wedge d\bar{z}_n \\
      &= P_n \,\norm{f(z)}_E^2 dx_1 \wedge dy_1 \wedge \cdots \wedge dx_n \wedge dx_n ,
   \end{split}
 \end{equation*}
where the constant $P_n = {(-2 \imath)}^n (-1)^{(n (n-1) / 2)}$ as before. Consequently, $\alpha \wedge \alpha^\star$ is invariant under local coordinate changes. In summary,

\thmcall{lemma}\label{lemma7.2}
Let $\Omega$, $\alpha$, $f$, and $P_n$ be as above. Then $\alpha \wedge \alpha^\star / P_n$ is a positive $2n$ dimensional, coordinate invariant volume measure on the under lying real manifold of $\Omega$, and
\begin{equation*}
  \alpha \wedge \alpha^* / P_n = \norm{f(z)}^2 \,dx_1 dy_1 \wedge \dots \wedge dx_n \wedge dy_n .
\end{equation*}
\exitthmcall{lemma}

\thmcall{definition}\label{def7.3}
  Let $\Omega$ be a connected, paracompact, analytic manifold of dimension $n \geq 1$ with basepoint $p \in \Omega$..  Let $E$ be a complex Hilbert space.  Let $k$ be an admissible weight function on $\Omega$.  Then the $E$ valued weighted Bergman space of differentials on $\Omega$, $\bbh = A^2(E, \Omega, k)$, is the space of all $E$ valued analytic $n$ forms $\alpha$ on $\Omega$ such that the weighted $\bbh$ norm
  \begin{equation*}
    \norm{\alpha}^2 =  \int_\Omega k \,\alpha \wedge \alpha^\star /P_n < \infty .
  \end{equation*}
\exitthmcall{definition}

Let $R = H^\infty(\Omega)$,  the algebra of bounded complex valued analytic functions on $\Omega$, or, in the case where $\Omega \subseteq \bbc^n$, let $R = \poly{\bbc}{z_1, \dots, z_n}$.  Define the order functions $\ord(r)$ and $\ord(\alpha)$ in the usual way, as the order of the zero of $r \in R$ or order of the zero of the analytic function $f$ at $p$, respectively. Here, $f$ is the analytic function appearing in the representation of $\alpha \in \bbh$ in local coordinates.

\

The definition of an analytic Hilbert module of $n$ forms is entirely analogous to definition \ref{def3.1} of an analytic module of functions:

\thmcall{definition}\label{def7.4}
$\bbh$ is an analytic Hilbert module of analytic $n$ forms on $\Omega$ over the analytic algebra $R$ if \,$\bbh$ is a left Hilbert module over $R$ under pointwise operations, and if each sequence $\alpha_n,\, n = 1, 2, 3, \dots$ in $\bbh$ converging to $\alpha$ in $\bbh$ norm also converges to $\alpha$ uniformly on compact subsets of $\Omega$. In other words, if the Hilbert space topology on $\bbh$ is finer than the topology of compact convergence on $\bbh$.
\exitthmcall{definition}

Of course, $\alpha_n$ converges to $\alpha$ uniformly on compact subsets of $\Omega$ if the analytic functions $f_n$ converge uniformly on compact subsets to $f$.  Here, $f_n$ and $f$ are the analytic functions appearing in the representations of $\alpha_n$ and $\alpha$ in local coordinates.

\thmcall{theorem}\label{thm7.5}
    $\bbh = A^2(E, \Omega, k \dmu)$ is an analytic Hilbert module of $E$ valued analytic differential $n$ forms over the analytic algebra $R$. Thus our abstract Beurling's theorem, \ref{thm2.5} and its corollary apply to $\bbh$, and so completely characterize the closed $R$ invariant subspaces of $\bbh$.
\exitthmcall{theorem}

\begin{proof}
The proof that $\bbh$ is an analytic Hilbert module of analytic differential $n$ forms is nearly identical to the proof of theorem \ref{thm5.2}, that $A^2(E, \Omega, k \dmu)$ is a an analytic Hilbert module of analytic functions.  And the proof that an analytic Hilbert module of analytic differential $n$ forms is a valuation Hilbert module is nearly identical to the proof of theorem \ref{thm3.2}, that an analytic Hilbert module of analytic functions is a valuation Hilbert module.
\end{proof}

\section{The Classical One Variable Beurling's Theorem}\label{sec8}

   Surely, it should be possible to derive the classical Beurling's theorem for $ H^2 $ of the unit disk from an abstract Beurling's theorem worthy of the name.  In this section, we are going to do just that.  In what follows, $ D $ shall be the unit disk, $ T $ the unit circle, $ \bbh = H^2 ( \bbc, D ) $, and $ V \neq \{ 0 \} $ shall be a closed subspace of $ \bbh $.
   
\thmcall{theorem}\label{thm8.1}
    Let $ \bbh  = H^2( \bbc, D ) $, and suppose the near homogeneous decomposition of $ V $ is near inner and satisfies the full projection property.  Then $ V $ is an $ R $ invariant subspace\! by theorem \ref{thm2.5} and its corollary.  Furthermore, there is a (complex valued) inner function $ I $ such that $ V = I \cdot H^2( \bbc, D ) $.
    
    Conversely, let $ I $ be an inner function, and suppose $ V = I \cdot \bbh $.  Then $ V $ is a closed $ R $ invariant subspace of $ \bbh $.
\exitthmcall{theorem}

\begin{proof}
As usual, let $ H_1 \plusperp H_2 \plusperp H_3 \plusperp \cdots $ be the near homogeneous decomposition of  $ \bbh = H^2( \bbc, D ) $.  The key to the proof is that each $ H_k $ is 1-dimensional, and in fact = $ \LH ( z^k) $, where $ \LH $ denotes the linear hull.  Recall from the notation introduced in paper I that $ P^\bbh_k $ is orthogonal projection from $ \bbh $ onto $ H_k $ (paper I, equations 5.5).  Let $ L_{ W_k } = P^\bbh_k $ restricted to $W_k$.  By paper I, proposition 5.17, $ L_{ W_k } $ is an invertible bounded linear transformation provided  $ W_k \neq \{ 0 \} $.

\

We are going to use the familiar isometric equivalence of $ \bbh $ with $ H^2 (\bbc, T )$, the space of square integrable complex valued functions defined a.e.\@ on $ T $, whose negative Fourier coefficients vanish.  Of course, as the reader will doubtless recall, this equivalence is defined by considering the a.e.\@ defined non-tangential boundary values of functions in $ \bbh $. We shall commit the gloss of identifying $ f \in \bbh $ with the corresponding a.e.\@ defined function $ f \in H^2 ( \bbc, T )$, and we shall let $ \bbh $ denote both Hilbert spaces interchangeably..

\

Now let $ k $ be the first index such that $ W_k \neq \{ 0 \} $.  Let $ I = I_k = L_{ W_k }^{-1} ( z^k ) $.  Since $ I \in W_k $,  $\ord( I ) = k $.  For each $ l  >  0 $, $ \ord( z^l I )  = k + l $, so $ z^I I \in V_{ k + l } $ .  But $ W_k \perp V_{ k + l } $, so $ I \perp z^l I $.  Thus  
\begin{equation*}
   0 = \inprod{ z^l I }{ I } = \int_0^{ 2 \pi } e^{ \imath l  \theta } I ( \eitheta ) \overline{ I ( \eitheta ) } \,\dtheta = \int_0^{ 2 \pi } e^{ \imath l \theta } \abs{ I ( \eitheta ) }^2 \,\dtheta ,
\end{equation*}
that is the $ l^\text{th} $ Fourier coefficient of the function $ \abs{ I }^2 $ is $ 0 $ for $ l = 1, 2, 3, \dots $.  Since $ \abs{ I }^2 $ is positive, the negative Fourier coefficients are also $ 0 $.  This $ I ( \eitheta ) $ is constant a.e.\@ on the unit circle $ T $.  By multiplying by a constant if necessary, $ \abs{ I ( \eitheta ) } = 1 $ a.e.\@ on $ T $.  Thus $ I $ is an inner function.

\

Because multiplication by $ z $ and $ I $ are isometries,  $ I \cdot \bbh $ is a closed $ R $ invariant subspace of $ H^2 ( \bbc, D ) $.  But we still have to prove that $ I \cdot \bbh $ is all of $ V $.  To that end, let $ j > k $ and consider $ W_j $.  As above, $ W_j $ is 1-dimensional, and $ W_j = \LH ( J ) $, where $ J $ is an inner function.  The question is, is $ J = z^{ j - k } I $?  The answer is yes for the following reasons:  First,  $ \ord (z^{ j - k } I )  = j $, so $ z^{ j - k } I \in  V_j $.   Second, the near homogeneous decomposition of $ V $ is near inner, so $ z^{  j - k + l } I \perp V_{ j  + l } $ for $ l = 1, 2, 3, \dots $  Thus $ z^{ j - k } \in W_j $.  Because $ W_j $ is 1-dimensional, by multiplying by a constant of modulus 1 if necessary, $ z^{ j - k } I =  J $.  Thus $ \LH ( z^m I ) = W_m $ for $m = k, k + 1, k + 2, \dots$, and so
\begin{equation*}
   \begin{split}
      I \cdot \bbh &= \LH ( z^k I ) \plusperp \LH ( z^{ k + 1 } I )\plusperp \LH ( z^{ k + 2 } I ) \plusperp \cdots \\
                         & = W_k  \plusperp W_{ k + 1 } \plusperp W_{ k + 2 } \cdots = V ,
   \end{split}
\end{equation*}
where the first equality is obvious, and the last equality follows from the fact that the near homogeneous decomposition converges to all of $ V $ by paper I, proposition 5.13.

\

   Conversely, let $ I $ be an inner function and $ V = I \cdot \bbh$.  Then it is clear from the fact that the map $ f  \mapsto I \cdot f $, $ f \in \bbh $, is an isometry that $ V $ is a closed, $ R $ invariant subspace.
\end{proof}

   Peter Lax generalized Beurling's theorem to vector valued $H^2$  of the half plane \cite{Lax1}. Henry Helson \cite{Hel1}, and Kenneth Hoffman \cite{Hoff1} translated Lax's theorem to the unit disk, that is to $H^2(E, D^1)$.  Our abstract Beurling's theorem adds nothing to Hoffman's proof of the vector valued case, so we shall not replicate it here.

\section{Extreme examples.}\label{sec9}
    Throughout this short section, $\bbh$ will be an analytic Hilbert module over an analytic algebra $R$.  It can happen that $\bbh$ and $R$ consist only of constants.  It can also happen that $\bbh$ contains non-constant analytic functions but $R$ consists only of constants.  We want to give examples of both situations and examine what our abstract Beurling's theorem, theorem \ref{thm2.5}, says in such extreme cases.

    There are easy artificial examples of these phenomena.  We can let $\bbh$ be any analytic Hilbert module containing constants and then change $\bbh$ to be those constant functions and $R$ to $\bbc$.  Or, we can let $\bbh$ be any analytic Hilbert module containing non-constant analytic functions but change $R$ to $\bbc$.  However, these examples are extremely artificial.  More telling are examples involving $H^2$ spaces for $\bbh$ and $H^\infty$ for $R$.  In the examples that follow, we shall return to using our previous notation involving $\Omega$ for domains in $\bbc^n$, $p$ for distinguished base points, etc. 

\thmcall{example}\label{example9.1}
$\Omega$ be the complex plane, $\Omega = \bbc$. Let $\bbh = H^2(E, \Omega)$, $R = H^\infty(\bbc, \Omega)$, $p = 0$.  In this case, $\bbh$ and $R$ consist entirely of constants.
\exitthmcall{example}

    It is clear why every function in $R$ is constant: every bounded entire function is constant by Liouville's theorem.  That $\bbh$ consists entirely of constants requires a brief argument:

    Let $f \in H^2(E, \Omega)$.  Then the subharmonic function $s_f \colon z \mapsto \norm{f(z)}_E^2$ has a harmonic majorant $h_f$.  We shall show that this implies each component of $f$ is constant.  So let $e \in E$ with $\norm{e}_E = 1$.  Then $\abs{\inprod{f(z)}{e}_E} \leq \norm{f(z)}_E$, and so the subharmonic function $s_e \colon z \mapsto \abs{\inprod{f(z)}{e}}_E^2 \leq s_f$.  Thus we have $s_{e} \leq s_f \leq h_f$, so $h_f - s_e$ is a positive subharmonic function on $\Omega$ and so is constant.  Hence $s_e$ is constant, so the analytic function $z \mapsto \inprod{f(z)}{e}_E$ has constant modulus.  By the maximum principle (or the open mapping theorem) for analytic functions, $z \mapsto \inprod{f(z)}{e}_E$ is constant.  Since $e$ was an arbitrary element of $E$ of norm $1$, each component of $f$ is constant, so $f$ itself is constant.\qed

\

    In the case of example \ref{example9.1}, our abstract Beurling's theorem still applies but says virtually nothing.  A closed $R$ submodule of $\bbh$ is simply a closed subspace of $E$, and vice versa.  If $V$ is a closed $R$ submodule, the valuation subspace series for $V$ is $V = V_0 \supseteq \{0\} \supseteq \{0\} \supseteq \cdots$, that is $V_1 = \{0\}$, $V_2 = \{0\}, \dots$ because, if a constant function vanishes at $0$, it is $0$.  The near homogeneous decomposition of $V$ is $V = V \plusperp \{0\} \plusperp \{0\} \plusperp \cdots$, that is $W_0 = V$, $W_1 = \{0\}$, $W_2 = \{0\} \dots\; $.  The reader can easily check that the near homogeneous decomposition of $V$ is near inner (because $R_1 = \{0\}$) and satisfies the full projection property (vacuously).  But these properties say nothing at all (beyond, perhaps, that $V$ is a closed subspace) in this extreme case.

\

    Now we shall temporarily let $p$ be an exponent rather than a basepoint.  Maurice Heins has given elegant examples of Riemann surfaces $\Omega$ such that $H^\infty(\bbc, \Omega)$ reduces to constant functions but $H^2(\bbc, \Omega)$ contains non-constant functions \cite[section III]{H2}.  In fact, he gave a finer classification than this: For each $p$, with $1 \leq p < \infty$, there are Riemann surfaces $\Omega$ such that $H^p(\bbc, \Omega)$ contains non-constant functions but $H^q(\bbc, \Omega)$ contains only constants for $p < q \leq \infty$ \cite{H2}.  (Here $q$ is simply an exponent, $1/p + 1/q$ is not necessarily equal to $1$.)  The $p = 2$, $q = \infty$ case will suffice for our purposes. 

\

    We now return to the convention that $p$ is a distinguished basepoint, not an exponent.

\thmcall{example}\label{example9.2}
Let $\Omega$ be a Riemann surface such that $H^2(\bbc, \Omega)$ contains non-constant analytic functions, but $H^\infty(\bbc, \Omega)$ consists entirely of constants.  Let $\bbh = H^2(\bbc, \Omega)$, $R = H^\infty(\bbc, \Omega)$, and $p \in \Omega$.
\exitthmcall{example}

In this case, a closed $R$ submodule of $\bbh$ is simply a closed subspace of $\bbh$, and vice versa.  Once again, our abstract Beurling's theorem applies but says virtually nothing.

\section{Open Questions.}\label{sec10}

First and foremost, can we deduce the Aleman, Richter, and Sunberg (ARS) characterization of the closed invariant subspaces of $A^2 (\bbc, D^1, \dmu)$ from our general abstract Beurling's theorem? More generally,

\thmcall{question}\label{quest10.1}
    For what dimensions $n$, for what domains or even complex analytic manifolds $\Omega$, and for what weighted $A^2$ spaces $A^2(\bbc, \Omega, k\dmu)$, can we deduce the ARS characterization of closed invariant subspaces, $V = \cl R \cdot M$, where $R = \poly{\bbc}{z_1, \ldots z_n}$ or $H^\infty(\Omega)$, and $M = V \orthcomp R \cdot V$, from our abstract Beurling's theorem?
\exitthmcall{question}

One of the simplest versions of question \ref{quest10.1} asks for which admissible weight functions $k$ can we deduce the ARS characterization of the closed invariant subspaces of $A^2(\bbc, D^1(0, 1),  k\dmu)$ from our abstract Beurling's theorem?  Shimorin \cite{Sh2} has conjectured that the ARS characterization of the closed invariant subspaces of the weighted Bergman space $A^2(\bbc, D^1, k\dmu)$ fails for the weight functions $k(z) = ( 1 + \alpha ) ( 1 - \abs{z}^2 )^\alpha$, with $\alpha > 4$. If the Shimorin conjecture is correct, the answer to this question is \emph{not for all}.\\

Of course, we can ask question \ref{quest10.1} even more generally:

\thmcall{question}\label{quest10.2}
    For what analytic algebras $R$ and analytic Hilbert modules $\bbh$, including vector valued ones, can we deduce the ARS characterization of closed invariant subspaces?
\exitthmcall{question}

\thmcall{question}\label{quest10.3}
    Halmos's paper also characterized the closed invariant subspaces of $L^2$ of the unit circle. See also the treatments in the books by Helson \cite{Hel1} and Hoffman \cite{Hoff1}. How can we generalize our abstract Beurling's theorem to the case of $L^2$ spaces in several complex variables?
\exitthmcall{question}

\thmcall{question}\label{quest10.4}
    How can we characterize the closed invariant subspaces of $H^p(E, \Omega)$ and $A^p(E, \Omega,$ \,$k\dmu)$ for $1 \leq p <  \infty$, and of $H^\infty(\Omega)$ under the $\beta$ topology \cite{Bu1}, in several complex variables?  In $1$ complex variable, the classical Beurling's theorem generalizes to these $H^p$ spaces (\cite{Hel1}, p.\@ 25).
\exitthmcall{question}

\thmcall{question}\label{quest10.5}
    Recall that a weakly near inner decomposition is one which satisfies the $g = 0$ case of the definition of a near inner decomposition (paper I, definition 6.5).  For which valuation Hilbert modules are all the weakly near inner, near homogeneous decompositions near inner?  Is this true for all valuation Hilbert modules?
\exitthmcall{question}

\end{document}